\documentclass[12pt, reqno]{amsart}
\usepackage{amsfonts}
\usepackage{amssymb}
\usepackage{amsmath}

\newtheorem{theorem}{Theorem}
\theoremstyle{plain}

\newtheorem{corollary}[theorem]{Corollary}

\newtheorem{proposition}[theorem]{Proposition}

\begin{document}
\title[] {On the fermionic $p$-adic integral representation of Bernstein polynomials associated with Euler numbers and polynomials}
\author{T. Kim}
\address{Taekyun Kim. Division of General Education-Mathematics \\
Kwangwoon University, Seoul 139-701, Republic of Korea  \\}
\email{tkkim@kw.ac.kr}
\author{J. Choi}
\address{Jongsung Choi. Division of General Education-Mathematics \\
Kwangwoon University, Seoul 139-701, Republic of Korea  \\}
\email{jeschoi@kw.ac.kr}
\author{Y.-H. Kim}
\address{Young-Hee Kim. Division of General
Education-Mathematics\\
Kwangwoon University, Seoul 139-701, Republic of Korea  \\}
\email{yhkim@kw.ac.kr}
\author{C. S. Ryoo}
\address{Cheon Seoung Ryoo. Department of Mathematics \\
Hannam University, Daejeon 306-791, Republic of Korea \\}
\email{ryoocs@hnu.kr}
%\thanks{$*$ correponding author}
\maketitle

{\footnotesize {\bf Abstract} \hspace{1mm} {
The purpose of this paper is to give some properties of several Bernstein type polynomials to represent the fermionic $p$-adic integral on $\Bbb Z_p$ .
From these properties, we derive some interesting identities on the Euler numbers and polynomials.
}

\section{Introduction}

Throughout this paper, let $p$ be an add prime number. The symbol
$\Bbb Z_p$, $\Bbb Q_p$ and $\Bbb C_p$ denote the ring of $p$-adic
integers, the field of $p$-adic numbers and the field of $p$-adic
completion of the algebraic closure of $\Bbb Q_p$, respectively (see
[1-14]).

Let $\Bbb N$ be the set of natural numbers and $\Bbb Z_+ =\Bbb N
\cup \{0\}$. Let $\nu_p$ be the normalized exponential valuation of
$\Bbb C_p$ with $|p|_p=p^{-\nu_p (p)} =\frac{1}{p}$. Note that
$\displaystyle \Bbb Z_p =\{x \ |\ |x|_p \le 1\}=
\lim_{\overleftarrow N} \Bbb Z \big/ p^N \Bbb Z_p $.

When one talks of $q$-extension, $q$ is variously consiered as an
indeterminate, a complex number $q \in \Bbb C$ or $p$-adic number $q
\in \Bbb C_p$. If $q\in \Bbb C$, we normally assume $|q|<1$ and if
$q\in \Bbb C_p$, we always assume $|1-q|_p <1$.

We say that $f$ is uniformly differentiable function at a point $a\in \Bbb Z$ and write $f\in UD(\Bbb Z_p)$,
if the difference quotient $F_f (x,y) =\frac{f(x)-f(y)}{x-y}$ have a limit $f'(a)$ as $(x,y)\rightarrow (a,a)$.
For $f\in UD(\Bbb Z_p)$, the fermionic $p$-adic $q$-integral on $\Bbb Z_p$ is defined as
\begin{eqnarray}
I_{-q} (f)=\int_{\Bbb Z_p } f(x) d\mu_{-q} (x) = \lim_{N \rightarrow
\infty} \frac{1+q}{1+q^{p^N}} \sum_{x=0}^{p^N-1} f(x)(-q)^x, \quad
(\text{see [7]}).
\end{eqnarray}
In the special case $q=1$ in (1), the integral
\begin{eqnarray}
I_{-1} (f)=\int_{\Bbb Z_p } f(x) d\mu_{-1} (x),
\end{eqnarray}
is called the fermionic $p$-adic invariant integral on $\Bbb Z_p$
(see [12]). From (2), we note
\begin{eqnarray}
I_{-1} (f_1)=-I_{-1}(f) + 2f(0),
\end{eqnarray}
where $f_1 (x)=f(x+1)$.

Moreover, for $n \in \Bbb N$, let $f_n (x)=f(x+n)$. Then we note
that
\begin{eqnarray}
I_{-1} (f_n)=(-1)^n I_{-1}(f) + 2\sum_{l=0}^{n-1}(-1)^{n-1-l}f(l).
\end{eqnarray}

It is well known that the Euler polynomials are defined by
\begin{eqnarray}
\frac{2}{e^t +1}e^{xt} = \sum_{n=0}^{\infty} E_n (x) \frac{t}{n!},
\quad \text{(see [1-14])}.
\end{eqnarray}
In the special case, $x=0$, $E_n(0)=E_n$ are called the $n$-th Euler
numbers.

Let $f(x)=e^{tx}$. Then, by (3), (4) and (5), we see that
\begin{eqnarray}
\int_{\Bbb Z_p} e^{(x+y)t} d\mu_{-1} (y) = \frac{2}{e^t +1}e^{xt} =
\sum_{n=0}^{\infty} E_n (x) \frac{t^n}{n!}.
\end{eqnarray}

Let $C[0, 1]$ denote the set of continuous functions on $[0, 1]$.
For $f \in C[0, 1]$, Bernstein introduced the following well-known
linear positive operator in the field of real numbers $\Bbb R$ :
\begin{eqnarray}
\mathbb{B}_n(f : x)= \sum_{k=0}^{n} f(\frac{k}{n})\binom{n}{k} x^k
(1-x)^{n-k}=\sum_{k=0}^{n} f(\frac{k}{n})B_{k,n}(x),
\end{eqnarray}
where
$\binom{n}{k}=\frac{n(n-1)\cdots(n-k+1)}{k!}=\frac{n!}{k!(n-k)!}$
(see [1, 2, 5, 9, 10, 14]). Here, $\mathbb{B}_n (f:x)$ is called the
Bernstein operator of order $n$ for $f$.

For $k,n \in \Bbb Z_+$, the Bernstein polynomial of degree $n$ is
defined by
\begin{eqnarray}
B_{k,n}= \binom{n}{k}x^k (1-x)^{n-k}, \quad \text{for} \ x\in [0,1].
\end{eqnarray}
For example, $B_{0,1}(x)=1-x$, $B_{1,1}(x)=x$, $B_{0,2}(x)=(1-x)^2$,
$B_{1,2}(x)=2x - 2x^2$, $B_{2,2}(x)= x^2, \ldots$, and
$B_{k,n}(x)=0$ for $n<k$, $B_{k,n}(x)=B_{n-k,n}(1-x)$.

In this paper, we study the properties of Bernstein polynomials in
the $p$-adic number field. For $f \in UD(\Bbb Z_p)$, we give some
properties of several type Bernstein polynomials to represent the
fermionic $p$-adic invariant integral on $\Bbb Z_p$. From those
properties, we derive some interesting identities on the Euler
polynomials.

\section{Fermionic $p$-adic integral representation of Bernstein polynomials}

By (5) and (6), we see that
\begin{eqnarray}
\frac{2}{e^t +1}e^{(1-x)t} =\sum_{n=0}^{\infty} E_n(1-x)
\frac{t^n}{n!}.
\end{eqnarray}
We also have that
\begin{eqnarray}
\frac{2}{e^t +1}e^{(1-x)t} =\frac{2}{1+e^{-t}}e^{-xt}=\sum_{n=0}^{\infty} E_n(x)\frac{(-1)^n}{n!}t^n .
\end{eqnarray}
From (9) and (10), we note tht $E_n(1-x) = (-1)^n E_n(x)$. It is
easy to show that
\begin{eqnarray}
E_n(2)=2- \sum_{l=0}^{n}\binom{n}{l} E_l=2 + E_n, \quad \text{for}
\,\, n>0.
\end{eqnarray}
By (5), (6), (9), (10) and (11), we see that for $n>0$,
\begin{eqnarray*}
\int_{\Bbb Z_p} (1-x)^n d\mu_{-1}(x)
&=& (-1)^n \int_{\Bbb Z_p} (x-1)^n d\mu_{-1}(x) \notag \\
&=& \int_{\Bbb Z_p} (x+2)^n d\mu_{-1}(x)  \\
&=& 2+ \int_{\Bbb Z_p} x^n d\mu_{-1}(x).\notag
\end{eqnarray*}
Therefore, we obtain the following theorem.

\begin{theorem}
For $n \in \Bbb N$, we have
\begin{eqnarray*}
\int_{\Bbb Z_p} (1-x)^n d\mu_{-1}(x)
= 2+ \int_{\Bbb Z_p} x^n d\mu_{-1}(x).
\end{eqnarray*}
\end{theorem}

\medskip

Theorem 1 is important to derive our main result in this paper.

Taking the fermionic $p$-adic integral on $\Bbb Z_p$ for one
Bernstein polynomial in (8), we get

\begin{eqnarray*}
\int_{\Bbb Z_p} B_{k,n}(x) d\mu_{-1}(x)
&=& \int_{\Bbb Z_p} \binom{n}{k} x^k (1-x)^{n-k} d\mu_{-1} (x)\\
&=& \binom{n}{k} \sum_{j=0}^{n-k} \binom{n-k}{j} (-1)^{n-k-j} \int_{\Bbb Z_p} x^{n-j} d\mu_{-1}(x)  \\
&=& \binom{n}{k} \sum_{j=0}^{n-k} \binom{n-k}{j} (-1)^{n-k-j} E_{n-j} \\
&=& \binom{n}{k} \sum_{j=0}^{n-k} \binom{n-k}{j} (-1)^{j} E_{k+j}.  \\
\end{eqnarray*}
Therefore, we obtain the following proposition.

\begin{proposition}
For $k, n \in \Bbb Z_+$, we have
\begin{eqnarray*}
\int_{\Bbb Z_p} B_{k,n}(x) d\mu_{-1}(x)
= \binom{n}{k} \sum_{j=0}^{n-k} \binom{n-k}{j}(-1)^j E_{k+j}.
\end{eqnarray*}
\end{proposition}

\medskip

It is known that $B_{k,n}(x)=B_{n-k,n}(1-x)$. Thus, we have
\begin{eqnarray}
\int_{\Bbb Z_p} B_{k,n}(x) d\mu_{-1}(x) &=& \int_{\Bbb Z_p}
B_{n-k,n}(1-x) d\mu_{-1}(x) \\
&=& \binom{n}{n-k} \sum_{j=0}^{k} \binom{k}{j} (-1)^{k-j} \int_{\Bbb
Z_p} (1-x)^{n-j} d\mu_{-1}(x).\notag
\end{eqnarray}

By (12) and Theorem 1, we see that for $n>k$,
\begin{eqnarray}
\int_{\Bbb Z_p} B_{k,n}(x) d\mu_{-1}(x) &=& \binom{n}{k}
\sum_{j=0}^{k} \binom{k}{j} (-1)^{k-j} (2+ \int_{\Bbb Z_p} x^{n-j}
d\mu_{-1}(x))\notag \\
&=& \binom{n}{k} \sum_{j=0}^{k} \binom{k}{j} (-1)^{k-j} (2+
E_{n-j}) \\
 &=& \left\{
\begin{array}{ll} 2+E_n \ \ &\hbox{if}\ \ k=0,
\vspace{2mm}\\
\binom{n}{k}  \sum_{j=0}^{k} \binom{k}{j} (-1)^{k-j}  E_{n-j} \ \
&\hbox{if}\ \ k>0.
\end{array}\right.  \notag
\end{eqnarray}
From (13), we obtain the following theorem.

\begin{theorem}
For $n, k \in \Bbb Z_+$ with $n >k$, we have
\begin{eqnarray*}
\int_{\Bbb Z_p} B_{k, n}(x) d\mu_{-1}(x) =\left\{
\begin{array}{ll} 2+E_n \ \ &\hbox{if}\ \ k=0,
\vspace{2mm}\\
\binom{n}{k} \sum_{j=0}^{k} \binom{k}{j} (-1)^{k-j}  E_{n-j} \ \
&\hbox{if}\ \ k>0.
\end{array}\right.
\end{eqnarray*}
\end{theorem}

\medskip

By Proposition 2 and Theorem 3, we obtain the following corollary.

\begin{corollary}
For $n, k \in \Bbb Z_+$ with $n >k$, we have
\begin{eqnarray*}
\sum_{j=0}^{n-k} \binom{n-k}{j} (-1)^{j}  E_{k+j} =\left\{
\begin{array}{ll} 2+E_n \ \ &\hbox{if}\ \ k=0,
\vspace{2mm}\\
\sum_{j=0}^{k} \binom{k}{j} (-1)^{k-j}  E_{n-j} \ \ &\hbox{if}\ \
k>0.
\end{array}\right.
\end{eqnarray*}
\end{corollary}

\medskip

For $m, n, k \in \Bbb Z_+$ with $m+n>2k$, fermionic $p$-adic
invariant integral for multiplication of two Berstein polynomials on
$\Bbb Z_p$ can be given by the following relation.
\begin{eqnarray*}
& &\int_{\Bbb Z_p} B_{k, n}(x)B_{k, m}(x) d\mu_{-1}(x)
\\ & & \qquad =\int_{\Bbb Z_p} \binom{n}{k}x^k (1-x)^{n-k} \binom{m}{k}x^k
(1-x)^{m-k} d \mu_{-1} (x)\\ & & \qquad = \binom{n}{k}\binom{m}{k}
\int_{\Bbb Z_p}x^{2k} (1-x)^{n+m-2k} d \mu_{-1} (x) \\ & & \qquad =
\binom{n}{k}\binom{m}{k} \sum_{j=0}^{2k}
\binom{2k}{j}(-1)^{j+2k}\int_{\Bbb Z_p} (1-x)^{n+m-j}d \mu_{-1} (x)
\\ & & \qquad =\binom{n}{k}\binom{m}{k} \sum_{j=0}^{2k}
\binom{2k}{j}(-1)^{j+2k}(2+\int_{\Bbb Z_p}
x^{n+m-j}d \mu_{-1} (x))\\
& & \qquad =\left\{
\begin{array}{ll} 2+E_{n+m} \ \ &\hbox{if}\ \ k=0,
\vspace{2mm}\\
\binom{n}{k}\binom{m}{k}\sum_{j=0}^{2k} \binom{2k}{j} (-1)^{j+2k}
E_{n+m-j} \ \ &\hbox{if}\ \ k>0.
\end{array}\right.
\end{eqnarray*}
Therefore, we obtain the following theorem.

\begin{theorem}
For $m, n, k \in \Bbb Z_+$ with $m+n>2k$, we have
\begin{eqnarray*}
& &\int_{\Bbb Z_p} B_{k, n}(x)B_{k, m}(x) d\mu_{-1}(x)
\\& & \qquad =\left\{
\begin{array}{ll} 2+E_{n+m}\ \ &\hbox{if}\ \ k=0,
\vspace{2mm}\\
\binom{n}{k}\binom{m}{k}\sum_{j=0}^{2k} \binom{2k}{j} (-1)^{j+2k}
E_{n+m-j} \ \ &\hbox{if}\ \ k>0.
\end{array}\right.
\end{eqnarray*}
\end{theorem}

\medskip

For $m, n, k \in \Bbb Z_+$, we have
\begin{eqnarray*}
& &\int_{\Bbb Z_p} B_{k, n}(x)B_{k, m}(x) d\mu_{-1}(x)
\\ & & \qquad = \binom{n}{k}\binom{m}{k}
\int_{\Bbb Z_p}x^{2k} (1-x)^{n+m-2k} d \mu_{-1} (x)
\\ & & \qquad =\binom{n}{k}\binom{m}{k} \sum_{j=0}^{n+m-2k}
\binom{n+m-2k}{j}(-1)^{j}\int_{\Bbb Z_p}
x^{j+2k}d \mu_{-1} (x)\\
& & \qquad =\binom{n}{k}\binom{m}{k} \sum_{j=0}^{n+m-2k}
\binom{n+m-2k}{j}(-1)^{j} E_{j+2k}.
\end{eqnarray*}
Thus, we obtain the following proposition.

\begin{proposition}
For $m, n, k \in \Bbb Z_+$, we have
\begin{eqnarray*}
\int_{\Bbb Z_p} B_{k, n}(x)B_{k, m}(x) d\mu_{-1}(x)
=\binom{n}{k}\binom{m}{k} \sum_{j=0}^{n+m-2k}
\binom{n+m-2k}{j}(-1)^{j} E_{j+2k}.
\end{eqnarray*}
\end{proposition}

\medskip

By Theorem 5 and Proposition 6, we obtain the following corollary.

\begin{corollary}
For $m, n, k \in \Bbb Z_+$ with $m+n>2k$, we have
\begin{eqnarray*}
& &\sum_{j=0}^{n+m-2k} \binom{n+m-2k}{j}(-1)^{j} E_{j+2k}
\\& & \qquad =\left\{
\begin{array}{ll} 2+E_{n+m}\ \ &\hbox{if}\ \ k=0,
\vspace{2mm}\\
\sum_{j=0}^{2k} \binom{2k}{j} (-1)^{j+2k} E_{n+m-j} \ \ &\hbox{if}\
\ k>0.
\end{array}\right.
\end{eqnarray*}
\end{corollary}

\medskip

In the same manner, multiplication of three Bernstein polynomials
can be given by the following relation :
\begin{eqnarray}
& &\int_{\Bbb Z_p} B_{k, n}(x)B_{k, m}(x)B_{k, s}(x) d\mu_{-1}(x)
\notag
\\ & &  = \binom{n}{k}\binom{m}{k}\binom{s}{k}
\sum_{j=0}^{n+m+s-3k} \binom{n+m+s-3k}{j}(-1)^{j}\int_{\Bbb Z_p}
x^{j+3k}d \mu_{-1} (x) \qquad \\
& &  =\binom{n}{k}\binom{m}{k}\binom{s}{k} \sum_{j=0}^{n+m+s-3k}
\binom{n+m+s-3k}{j}(-1)^{j} E_{j+3k}, \notag
\end{eqnarray}
where $m, n, s, k \in \Bbb Z_+$ with $m+n+s>3k$.

For $m, n, s, k \in \Bbb Z_+$ with $m+n+s>3k$, by the symmetry of
Bernstein polynomals, we see that
\begin{eqnarray*}
& &\int_{\Bbb Z_p} B_{k, n}(x)B_{k, m}(x)B_{k, s}(x) d\mu_{-1}(x)
\\ & & \qquad = \binom{n}{k}\binom{m}{k}\binom{s}{k}
\sum_{j=0}^{3k} \binom{3k}{j}(-1)^{3k-j}\int_{\Bbb Z_p}
(1-x)^{n+m+s-j}d \mu_{-1} (x)\\
& & \qquad =\binom{n}{k}\binom{m}{k}\binom{s}{k} \sum_{j=0}^{3k}
\binom{3k}{j}(-1)^{3k-j}(2+\int_{\Bbb Z_p} x^{n+m+s-j}\mu_{-1} (x))
\\& & \qquad =\left\{
\begin{array}{ll} 2+E_{n+m+s}\ \ &\hbox{if}\ \ k=0,
\vspace{2mm}\\
\binom{n}{k}\binom{m}{k}\binom{s}{k} \sum_{j=0}^{3k}
\binom{3k}{j}(-1)^{3k-j} E_{n+m+s-j} \ \ &\hbox{if}\ \ k>0.
\end{array}\right.
\end{eqnarray*}
Therefore, we obtain the following theorem.

\begin{theorem}
For $m, n, s, k \in \Bbb Z_+$ with $m+n+s>3k$, we have
\begin{eqnarray*}
& &\int_{\Bbb Z_p} B_{k, n}(x)B_{k, m}(x)B_{k, s}(x) d\mu_{-1}(x)
\\& & \qquad =\left\{
\begin{array}{ll} 2+E_{n+m+s}\ \ &\hbox{if}\ \ k=0,
\vspace{2mm}\\
\binom{n}{k}\binom{m}{k}\binom{s}{k} \sum_{j=0}^{3k}
\binom{3k}{j}(-1)^{3k-j} E_{n+m+s-j} \ \ &\hbox{if}\ \ k>0.
\end{array}\right.
\end{eqnarray*}
\end{theorem}

\medskip

By (14) and Theorem 8, we obtain the following corollary.

\begin{corollary}
For $m, n, s, k \in \Bbb Z_+$ with $m+n+s>3k$, we have
\begin{eqnarray*}
& &\sum_{j=0}^{n+m+s-3k} \binom{n+m+s-3k}{j}(-1)^{j} E_{j+3k}
\\& & \qquad =\left\{
\begin{array}{ll} 2+E_{n+m+s}\ \ &\hbox{if}\ \ k=0,
\vspace{2mm}\\
\sum_{j=0}^{3k} \binom{3k}{j} (-1)^{3k-j} E_{n+m+s-j} \ \
&\hbox{if}\ \ k>0.
\end{array}\right.
\end{eqnarray*}
\end{corollary}

Using the above theorems and mathematical induction, we obtain the
following theorem.

\begin{theorem} Let $s \in \Bbb N$.
For $n_1, n_2, \ldots, n_s, k \in \Bbb Z_+$ with $n_1+n_2+\cdots+
n_s>sk$, the multiplication of the sequence of Bernstein polynomials
$B_{k, n_1} (x), \ldots, B_{k, n_s}(x)$ with different degrees under
fermionic $p$-adic invariant integral on $\Bbb Z_p$ can be given as
\begin{eqnarray*}
& &\int_{\Bbb Z_p} \left(\prod_{i=1}^{s} B_{k, n_i}(x) \right)
d\mu_{-1}(x)
\\& & \qquad =\left\{
\begin{array}{ll} 2+E_{n_1+n_2+\cdots+
n_s}\ \ &\hbox{if}\ \ k=0,
\vspace{2mm}\\
\left(\prod_{i=1}^{s} \binom{n_i}{k} \right) \sum_{j=0}^{sk}
\binom{s k}{j}(-1)^{sk-j} E_{n_1+n_2+\cdots+ n_s-j} \ \ &\hbox{if}\
\ k>0.
\end{array}\right.
\end{eqnarray*}
\end{theorem}

\medskip

We also easily see that
\begin{eqnarray}
& &\int_{\Bbb Z_p} \left(\prod_{i=1}^{s} B_{k, n_i}(x) \right)
d\mu_{-1}(x) \\ & & \qquad =\left(\prod_{i=1}^{s} \binom{n_i}{k}
\right) \sum_{j=0}^{n_1 + \cdots +n_s - sk} \binom{n_1 + \cdots +n_s
- s k}{j}(-1)^{j} E_{j+s k}. \notag
\end{eqnarray}

By Theorem 10 and (15), we obtain the following corollary.

\begin{corollary} Let $s \in \Bbb N$.
For $n_1, n_2, \ldots, n_s, k \in \Bbb Z_+$ with $n_1+n_2+\cdots+
n_s>sk$, we have
\begin{eqnarray*}
& &  \sum_{j=0}^{n_1 + \cdots +n_s - sk} \binom{n_1 + \cdots +n_s -
s k}{j}(-1)^{j} E_{j+s k}
\\& & \qquad =\left\{
\begin{array}{ll}  2+E_{n_1+n_2+\cdots+
n_s}\ \ &\hbox{if}\ \ k=0,
\vspace{2mm}\\
\sum_{j=0}^{sk} \binom{s k}{j}(-1)^{sk-j} E_{n_1+n_2+\cdots+ n_s-j}
\ \ &\hbox{if}\ \ k>0.
\end{array}\right.
\end{eqnarray*}
\end{corollary}

\medskip

Let $m_1, \ldots, m_s, n_1, \ldots, n_s, k \in \Bbb Z_+$ with $m_1
n_1+ \cdots+ m_s n_s>(m_1+\cdots+ m_s)k$. By the definition of
$B_{k, n_s}^{m_s}(x)$, we easily get
\begin{eqnarray*}
& &\int_{\Bbb Z_p} \left(\prod_{i=1}^{s} B_{k, n_i}^{m_i}(x) \right)
d\mu_{-1}(x)
\\& & =\left(\prod_{i=1}^{s} \binom{n_i}{k}^{m_i} \right) \sum_{j=0}^{k \sum_{i=1}^s m_i}
(-1)^{k\sum_{i=1}^s m_i-j} \int_{\Bbb Z_p} (1-x)^{\sum_{i=1}^{s} n_i
m_i -j} d \mu_{-1}(x)
\\& & = \left(\prod_{i=1}^{s} \binom{n_i}{k}^{m_i} \right) \sum_{j=0}^{k \sum_{i=1}^s m_i}
\binom{k \sum_{i=1}^s m_i}{j}(-1)^{k\sum_{i=1}^s
m_i-j}(2+E_{\sum_{i=1}^{s} n_i m_i -j})
\\& & =\left\{
\begin{array}{ll}  2+E_{m_1
n_1+ \cdots+ m_s n_s}\ \ &\hbox{if}\ \ k=0,
\vspace{2mm}\\
\left(\prod_{i=1}^{s} \binom{n_i}{k}^{m_i} \right) \sum_{j=0}^{k
\sum_{i=1}^s m_i} \binom{k \sum_{i=1}^s m_i}{j}(-1)^{k\sum_{i=1}^s
m_i-j}E_{\sum_{i=1}^{s} n_i m_i -j} \ \ &\hbox{if}\ \ k>0.
\end{array}\right.
\end{eqnarray*}
Therefore, we obtain the following theorem.

\begin{theorem} Let $s \in \Bbb N$.
For $m_1, \ldots, m_s, n_1, \ldots, n_s, k \in \Bbb Z_+$ with $m_1
n_1+ \cdots+ m_s n_s>(m_1+\cdots+ m_s)k$, we have
\begin{eqnarray*}
& &\int_{\Bbb Z_p} \left(\prod_{i=1}^{s} B_{k, n_i}^{m_i}(x) \right)
d\mu_{-1}(x)\\
\\& & =\left\{
\begin{array}{ll}  2+E_{m_1
n_1+ \cdots+ m_s n_s}\ \ &\hbox{if}\ \ k=0,
\vspace{2mm}\\
\left(\prod_{i=1}^{s} \binom{n_i}{k}^{m_i} \right) \sum_{j=0}^{k
\sum_{i=1}^s m_i} \binom{k \sum_{i=1}^s m_i}{j}(-1)^{k\sum_{i=1}^s
m_i-j}E_{\sum_{i=1}^{s} n_i m_i -j} \ \ &\hbox{if}\ \ k>0.
\end{array}\right.
\end{eqnarray*}
\end{theorem}

\medskip

By simple calculation, we easily get
\begin{eqnarray*}
& &\int_{\Bbb Z_p} \left(\prod_{i=1}^{s} B_{k, n_i}^{m_i}(x) \right)
d\mu_{-1}(x)
\\& &=\left(\prod_{i=1}^{s} \binom{n_i}{k}^{m_i} \right) \sum_{j=0}^{
\sum_{i=1}^s n_i m_i- k\sum_{i=1}^s m_i} \binom{\sum_{i=1}^s n_i
m_i- k\sum_{i=1}^s m_i}{j}(-1)^j E_{ k \sum_{i=1}^{s} m_i-j},
\end{eqnarray*}
where $m_1, \ldots, m_s, n_1, \ldots, n_s, k \in \Bbb Z_+$ for $s
\in \Bbb N$. By Theorem 12 and (16), we obtain the following
corollary.

\begin{corollary}
Let $s \in \Bbb N$. For $m_1, \ldots, m_s, n_1, \ldots, n_s, k \in
\Bbb Z_+$ with $m_1 n_1+ \cdots+ m_s n_s>(m_1+\cdots+ m_s)k$, we
have
\begin{eqnarray*}
& &   \sum_{j=0}^{ \sum_{i=1}^s n_i m_i- k\sum_{i=1}^s m_i}
\binom{\sum_{i=1}^s n_i m_i- k\sum_{i=1}^s m_i}{j}(-1)^j E_{ k
\sum_{i=1}^{s} m_i-j} \\
\\& & \qquad = \left\{
\begin{array}{ll}  2+E_{m_1
n_1+ \cdots+ m_s n_s}\ \ &\hbox{if}\ \ k=0,
\vspace{2mm}\\
\sum_{j=0}^{k \sum_{i=1}^s m_i} \binom{k \sum_{i=1}^s
m_i}{j}(-1)^{k\sum_{i=1}^s m_i-j}E_{\sum_{i=1}^{s} n_i m_i -j} \ \
&\hbox{if}\ \ k>0.
\end{array}\right.
\end{eqnarray*}
\end{corollary}

\medskip

The fermionic $p$-adic invariant integral of multiplication of
($n+1$) Bernstein polynomials, the $n$-th degree Bernstein
polynomials $B_{i,n}(x)$ with $i=0, 1, \ldots, n$ and with
multiplicity $m_0, m_1, \ldots, m_n$ on $\Bbb Z_p$, respectively,
can be given by
\begin{eqnarray*}
& &\int_{\Bbb Z_p} \left(\prod_{i=0}^{n} B_{i, n}^{m_i}(x) \right)
d\mu_{-1}(x) \\
\\& & =\left(\prod_{i=0}^{n} \binom{n}{i}^{m_i} \right)\int_{\Bbb
Z_p} x^{\sum_{i=1}^{n}i m_i}(1-x)^{n \sum_{i=0}^{n} m_i -
\sum_{i=1}^n i m_i} d \mu_{-1} (x) \\ \\& & = \frac{(\prod_{i=1}^{n}
\binom{n}{i}^{m_i})}{\binom{n\sum_{i=0}^{n} m_i}{\sum_{i=1}^{n}i
m_i}} \int_{\Bbb Z_p} B_{\sum_{i=1}^{n}i m_i, \, n\sum_{i=0}^{n}
m_i} (x) d \mu_{-1} (x),
\end{eqnarray*}
where $m_0, m_1, \ldots, m_n \in \Bbb Z_+$ with $n \in \Bbb Z_+$.

Assume that $nm_0 +n m_1+ \cdots+ n m_n >m_1+2m_2 +\cdots+ n m_n$.
Then we have
\begin{eqnarray*}
& &\int_{\Bbb Z_p} \left(\prod_{i=0}^{n} B_{i, n}^{m_i}(x) \right)
d\mu_{-1}(x) \\
\\& & =\left\{
\begin{array}{ll}  2+E_{nm_0 +n m_1+ \cdots+ n m_n }\qquad \qquad \qquad \qquad \qquad \qquad \qquad \text{if} \ \ \sum_{i=1}^{n}i
m_i=0,
\vspace{2mm}\\
\left(\prod_{i=0}^{n} \binom{n}{i}^{m_i} \right)\sum_{j=0}^{
\sum_{i=1}^n i m_i} \binom{\sum_{i=1}^n i m_i}{j}(-1)^{\sum_{i=1}^n
i m_i-j}E_{n\sum_{i=0}^{n} m_i -\sum_{i=1}^n i m_i} \vspace{2mm} \\
\qquad \qquad \qquad \qquad \qquad \qquad \qquad \qquad  \qquad
\qquad \qquad \qquad \text{if}\ \ \sum_{i=1}^{n}i m_i>0.
\end{array}\right.
\end{eqnarray*}
Therefore, we obtain the following theorem.

\begin{theorem} Let $n \in \Bbb Z_+$.\\
(I) For $m_0, m_1, \ldots, m_n \in \Bbb Z_+$ with $n\sum_{i=0}^{n}
m_i > \sum_{i=1}^n i m_i$, we have
\begin{eqnarray*}
& &\int_{\Bbb Z_p} \left(\prod_{i=0}^{n} B_{i, n}^{m_i}(x) \right)
d\mu_{-1}(x)
\\& & =\left\{
\begin{array}{ll}  2+E_{nm_0 +n m_1+ \cdots+ n m_n }\qquad \qquad \qquad \qquad \qquad \qquad \qquad \text{if} \ \ \sum_{i=1}^{n}i
m_i=0,
\vspace{2mm}\\
\left(\prod_{i=0}^{n} \binom{n}{i}^{m_i} \right)\sum_{j=0}^{
\sum_{i=1}^n i m_i} \binom{\sum_{i=1}^n i m_i}{j}(-1)^{\sum_{i=1}^n
i m_i-j}E_{n\sum_{i=0}^{n} m_i -\sum_{i=1}^n i m_i} \vspace{2mm} \\
\qquad \qquad \qquad \qquad \qquad \qquad \qquad \qquad  \qquad
\qquad \qquad \qquad \text{if}\ \ \sum_{i=1}^{n}i m_i>0.
\end{array}\right.
\end{eqnarray*}
(II) For $m_0, m_1, \ldots, m_n \in \Bbb Z_+$, we have
\begin{eqnarray*}
& &\int_{\Bbb Z_p} \left(\prod_{i=0}^{n} B_{i, n}^{m_i}(x) \right)
d\mu_{-1}(x)
\\& & =\left(\prod_{i=0}^{n} \binom{n}{i}^{m_i} \right)\sum_{j=0}^{
n\sum_{i=0}^{n} m_i -\sum_{i=1}^n i m_i} \binom{n\sum_{i=0}^{n} m_i
-\sum_{i=1}^n i m_i}{j}(-1)^{j}E_{\sum_{i=1}^n i m_i+j}.\\
\end{eqnarray*}
\end{theorem}

\medskip

By Theorem 14, we obtain the following corollary.
\begin{corollary}
For $n, m_0, m_1, \ldots, m_n \in \Bbb Z_+$ with $n\sum_{i=0}^{n}
m_i > \sum_{i=1}^n i m_i$, we have
\begin{eqnarray*}
& &   \sum_{j=0}^{ n\sum_{i=0}^{n} m_i -\sum_{i=1}^n i m_i}
\binom{n\sum_{i=0}^{n} m_i -\sum_{i=1}^n i
m_i}{j}(-1)^{j}E_{\sum_{i=1}^n i m_i+j}
\\ \\
& & =\left\{
\begin{array}{ll}  2+E_{nm_0 +n m_1+ \cdots+ n m_n }\qquad \qquad \qquad \qquad \qquad \qquad \qquad \text{if} \ \ \sum_{i=1}^{n}i
m_i=0,
\vspace{2mm}\\
\sum_{j=0}^{ \sum_{i=1}^n i m_i} \binom{\sum_{i=1}^n i
m_i}{j}(-1)^{\sum_{i=1}^n
i m_i-j}E_{n\sum_{i=0}^{n} m_i -\sum_{i=1}^n i m_i} \vspace{2mm} \\
\qquad \qquad \qquad \qquad \qquad \qquad \qquad \qquad  \qquad
\qquad \qquad \qquad \text{if}\ \ \sum_{i=1}^{n}i m_i>0.
\end{array}\right.
\end{eqnarray*}
\end{corollary}


\bigskip


\begin{thebibliography}{99}

\bibitem{1} M. Acikgoz, S. Araci, \textit{A study on the integral of the product of several type Berstein
polynomials}, IST Transaction of Applied Mathematics-Modelling and
Simulation, 2010.

\bibitem{2} S. Bernstein, \textit{D$\acute{e}$monstration du th$\acute{e}$or$\grave{e}$me de Weierstrass, fond$\acute{e}$e sur le calcul des probabilities},
Commun. Soc. Math. Kharkow \textbf{13} (1912), 1--2.


\bibitem{4} I. N. Cangul, V. Kurt, H. Ozden, Y. Simsek, \textit{On the higher-order $w$-$q$-Genocchi numbers},
Adv. Stud. Contemp. Math. \textbf{19} (2009), 39--57.

\bibitem{4} N. K. Govil, V. Gupta, \textit{Convergence of $q$-Meyer-K$\ddot{o}$nig-Zeller-Durrmeyer operators},
Adv. Stud. Contemp. Math. \textbf{19} (2009), 97--108.

\bibitem{5} V. Gupta, T. Kim, J. Choi, Y.-H. Kim, \textit{Generating function for $q$-Bernstein,  $q$-Meyer-K$\ddot{o}$nig-Zeller and $q$-Beta
basis}, Automation Computers Applied Mathematics \textbf{19} (2001),
7--11.

\bibitem{6} T. Kim, \textit{$q$-extension of the Euler formulae and trigonometric functions}, Russ. J. Math. Phys. \textbf{14} (2007), 275--278.

\bibitem{7} T. Kim, \textit{$q$-Volkenborn integration},  Russ. J. Math. Phys.   \textbf{9} (2002), 288--299.

\bibitem{8} T. Kim, \textit{$q$-Bernoulli numbers and polynomials associated with Gaussian binomial coefficients},
Russ. J. Math. Phys.   \textbf{15} (2008), 51--57.

\bibitem{9} T. Kim, J. Choi, Y.-H. Kim, \textit{Some identities on the $q$-Berstein polynomials, $q$-Stirling numbers and $q$-Bernoulli numbers},
Adv. Stud. Contemp. Math. \textbf{20} (2010), 335-341

\bibitem{10} T. Kim, L.-C. Jang, H. Yi, \textit{Note on the modified $q$-Berstein polynomials}, Discrete Dynamics in Nature and Society (in
press), (2010), arXiv 1005.4293.

\bibitem{11} T. Kim,  \textit{Note on the Euler $q$-zeta functions}, J. Number Theory  \textbf{129} (2009), 1798--1804.

\bibitem{12} T. Kim, \textit{Barnes-type multiple $q$-zeta functions and $q$-Euler polynomials},
J. Physics A : Math. Theor., \textbf{43} (2010), 255201, 11pp.

\bibitem{13} V. Kurt,\textit{A further symmetric relation on the analogue of the Apostol-Bernoulli and the analogue of the Apostol-Genocchi polynomials},
Appl. Math. Sci. \textbf{3} (2009), 53--56.

\bibitem{14} Y. Simsek, M. Acikgoz, \textit{A new generating function of
$q$-Berstein-type polynomials and their interpolation function},
Abstract and Applied Analysis \textbf{2010} (2010), Article ID
769095, 12 pages.





\end{thebibliography}
\end{document}